\documentclass[times]{nmeauthNoIJNME}

\usepackage{csml}

\usepackage{mathtools}

\usepackage[pdftex,colorlinks=true,allcolors=blue]{hyperref}

\usepackage{subfig}

\graphicspath{ {./figs/} }

\begin{document}

\runningheads{X Xiao et al.}{A non-iterative method for robustly computing intersections}

\title{A non-iterative method for robustly computing the intersections between a line and a curve or surface}

\author{Xiao Xiao\affil{1}, Laurent Bus\'e\affil{2} and Fehmi Cirak\affil{1}\corrauth}
\address{\affilnum{1}Department of Engineering, University of Cambridge, Trumpington Street, Cambridge CB2 1PZ, UK \break 
\affilnum{2}Universit\'e C\^ote d'Azur, Inria, 2004 route des Lucioles, 06902 Sophia Antipolis, France
}

\corraddr{Department of Engineering, University of Cambridge, Trumpington Street, Cambridge CB2 1PZ, UK. E-mail: f.cirak@eng.cam.ac.uk}

\begin{abstract}
The need to compute the intersections between a line and a high-order curve or surface arises in a large number of finite element applications. Such intersection problems are easy to formulate but hard to solve robustly. We introduce a non-iterative method for computing intersections by solving a matrix singular value decomposition (SVD) and an eigenvalue problem. That is, all intersection points and their parametric coordinates are determined in one-shot using only standard linear algebra techniques available in most software libraries. As a result, the introduced technique is far more robust than the widely used Newton-Raphson iteration or its variants. The maximum size of the considered matrices depends on the polynomial degree $q$ of the shape functions and is \mbox{$2q \times 3q$} for curves and \mbox{$6 q^2 \times 8 q^2$} for surfaces. The method has its origin in algebraic geometry and has here been considerably simplified with a view to widely used high-order finite elements. In addition, the method is derived from a purely linear algebra perspective without resorting to algebraic geometry terminology. A complete implementation is available from \href{http://bitbucket.org/nitro-project/}{http://bitbucket.org/nitro-project/}.
\end{abstract}


\maketitle

%\tableofcontents 

%
%--------------------------------------------------------------------------------
\section{Introduction \label{sec:intro}}
%--------------------------------------------------------------------------------
%
There has recently been an increased academic and industrial interest in high-order finite elements due to their efficiency advantages over classical low-order elements, see e.g.\cite{Hughes:2005aa,karniadakis2013spectral}. To achieve their full potential, high-order methods require the curved domain boundaries to be approximated with non-planar elements. The intersection between the curved elements and a given line is required in a wide range of applications, including contact~\cite{wriggers2004computational}, mesh generation~\cite{xie2013generation,turner2018curvilinear} and immersed finite elements~\cite{Ruberg:2011aa}. As known, intersection computations lead to an easy to formulate, but hard to solve nonlinear root-finding problem. The prevalent technique in computational mechanics for solving such problems is Newton-Raphson iteration, which is in general not very robust, especially when no good starting points are available or multiple intersections are present.

In computer-aided geometric design (CAD) and manufacturing (CAM) intersection computation is a recurring task and, to this end, a number of ingenious methods have been developed~\cite{patrikalakis2009shape}.  Especially promising are the non-iterative methods with an origin in algebraic geometry, which are for the most part unknown in computational mechanics. Algebraic geometry deals with systems of polynomial equations and geometric objects defined by them and provides the most rigorous framework for intersection computations~\cite{sederbergCN:2012}.  The specific technique considered in this paper is the implicitisation technique proposed by Bus\'e~\cite{laurent2014implicit}, which shares some commonalities with the method of moving lines/planes introduced by Sederberg and Chen~\cite{sederberg1995implicitization}. Different from the original work, in this paper we derive the method from a purely linear algebra viewpoint. To follow the presented derivations, it is sufficient to only know the notion of the algebraic degree of a curve or surface. The algebraic degree of a curve or surface is defined as its number of intersections with a line. Counting all intersections (real, complex, multiple ones and ones at infinity), the algebraic degree of a polynomial curve of degree~$q$ and the corresponding tensor-product surface are~$q$ and~$2q^2$, respectively.  As will become clear, adopting a linear algebra viewpoint has the added benefit that many linear algebra techniques, like preconditioning and various matrix decompositions, become readily available for  intersection computations.

In the following, we first discuss  the intersection between a planar Lagrange curve and a line and provide an easy to follow illustrative example. This simple case is sufficient to introduce and discuss the key aspects of the proposed non-iterative technique. Its  extension to the surface case is straightforward and is discussed last.

%!TEX root = /Users/lbuse/Inria/Article_Tex/nitro/paper/nonIterative.tex
%
%--------------------------------------------------------------------------------
\section{Intersection of lines with curves  \label{sec:curves}}
%--------------------------------------------------------------------------------
%

%
%--------------------------------------------------------------------------------
\subsection{Moving lines and intersections}
%--------------------------------------------------------------------------------
%
%
Let~\mbox{$\vec x (\theta) = \left ( x^1 (\theta), \, x^2 (\theta) \right)^\trans $} be a  planar parametric curve, with~$\vec x(\theta) \in \mathbb R^2$,  of degree~$q_x$ given either in Lagrange basis~$L_i(\theta)$ or power (monomial) basis~$P_j(\theta) $ with 
\begin{equation} \label{eq:curve}
	\vec  x (\theta )= \sum_{i=1}^{q_x+1}  L_i (\theta)   \vec x_i  =   \sum_{j=1}^{q_x+1} P_j(\theta)    \vec \alpha_j \,  ,
\end{equation}
where~$\vec x_i \in \mathbb R^2$ are the nodal coordinates and~$\vec \alpha_j \in \mathbb R^2$ are the coefficients in the power basis. As usual, the power basis~$P_j(\theta)$ contains the consecutive powers of~$\theta$ from~$0$ up to~$q_x$. The two basis are related by
\begin{equation}
	P_j(\theta) =  \sum_{i=1}^{q_x+1} L_i(\theta) P_j(\theta_i) \, ,
\end{equation}
where~$P_j(\theta_i)$ is the Vandermonde matrix and~$\theta_i$ is the parametric coordinate of the $i$-th Lagrange node.  Following a similar approach a curve given in any other polynomial basis can be re-expressed in the power basis.%

\begin{figure}
\centering
\subfloat[Two linear moving lines ($q_g = 1$) defining a quadratic curve ($q_x=2$)]
	{
		\includegraphics[scale = 0.9]{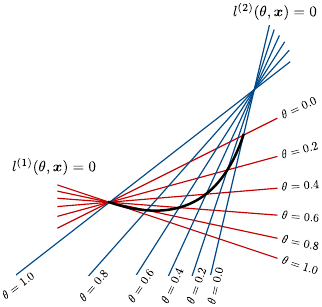}
		\label{fig:movingLines2}
	}
\hspace{1.5cm}
\subfloat[Four of the five cubic moving lines ($q_g=3$) defining a cubic curve ($q_x=3$) at $\theta \in \{ 0.4, \, 0.8\}$] {
		\includegraphics[scale = 0.9]{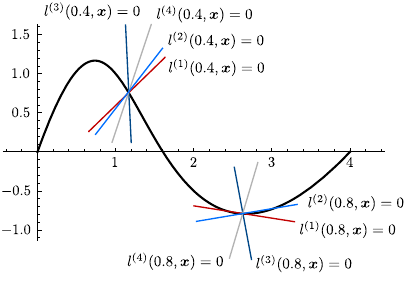}
		\label{fig:curveMovingLines}
	}
\caption{Moving lines and their intersections.}
\end{figure}

To define a point~$\vec x( \theta)$ on the curve as the intersection of several moving lines, or pencils of lines,  consider the line
\begin{equation}~\label{eq:lines}
	 l (\theta, \, \vec x) =
	\begin{pmatrix}
	\vec x  \\ 1
	\end{pmatrix} 
	 \cdot \vec g (\theta) =  x^1 g^1(\theta) +  x^2 g^2(\theta) +    g^3(\theta)=  0 \,  ,
\end{equation}
where $\vec x \in \mathbb R^2$ and ~$\vec g (\theta) = \left ( g^1(\theta), \, g^2 (\theta), \, g^3(\theta)
\right)^\trans$ is an auxiliary vector collecting the parameters of the line. For a fixed~$\theta$ equation~\eqref{eq:lines} describes a line and the line moves with the parameter~$\theta$, see Figure~\ref{fig:movingLines2} .   The three parameters of the line are assumed to be polynomial functions given by
\begin{equation}~\label{eq:aux}
	\vec g (\theta) = \sum_{l=1}^{q_g+1} \widetilde{P}_l (\theta) \vec g_l \, .
\end{equation}
The degree~$q_g$ of the power basis~$\widetilde P_l (\theta)$ has to be chosen sufficiently high in order to be able to compute all the intersection points (real, complex, multiple ones and ones at infinity). The number of intersection points is equivalent to the algebraic degree of the curve~$\vec x(\theta)$. A curve~$\vec x(\theta)$ of degree~$q_x$ has~$q_x$ intersection points with a line. As will become clear, the number of intersections implies a constraint on the minimum possible degree~$q_g$ for~$\widetilde P_l(\theta)$. 

Next, the aim is to find  several moving lines of the form~\eqref{eq:lines}, or more specifically their coefficients~$\vec g_l$ in~\eqref{eq:aux},  with a common intersection point which is a point on the curve~$\vec x(\theta)$. It is required that each line satisfies at the common intersection point
 \begin{equation}~\label{eq:linesDisc}
 	 l (\theta, \, \vec x(\theta) ) =
	\begin{pmatrix}
	\vec x (\theta) \\ 1
	\end{pmatrix} 
	 \cdot \vec g (\theta) =  
	  \begin{pmatrix}
	\vec x (\theta) \\ 1
	\end{pmatrix}  
	\cdot 
	\left (  \sum_{l=1}^{q_g+1}  
	\widetilde{P}_l  (\theta) \vec g_l	  \right )
	 =  0 \, .
\end{equation}
 After introducing the definition of the curve~\eqref{eq:curve} this yields 
\begin{equation}~\label{eq:linesDisc2}
	    \sum_{l=1}^{q_g+1} \left (   \sum_{j=1}^{q_x+1}  \begin{pmatrix} \vec \alpha_j   \\ 1 \end{pmatrix} P_j(\theta)   \widetilde{P}_l (\theta) \right)  \vec g_l  = 0 \, .
\end{equation}
The bracketed term can be expressed in a new power basis~$\widehat P_k$ of dimension~$q_x+q_g+1$ with 
\begin{equation}
	 \sum_{k=1}^{q_x+q_q+1}   \widehat P_k  (\theta) C_{kl} =   \sum_{j=1}^{q_x+1}  \begin{pmatrix}   \vec \alpha_j   \\ 1 \end{pmatrix}  P_j(\theta)   \widetilde{P}_l (\theta) \, , 
\end{equation}
where the matrix components~$C_{kl}$  contain the known coefficients~$\vec \alpha_j$. Equation~\eqref{eq:linesDisc2} can now be rewritten as
\begin{equation}~\label{eq:moveLineDesc}
	   \sum_{l=1}^{3(q_g+1)}  \left ( \sum_{k=1}^{q_x+q_g+1}  \widehat P_k (\theta)  C_{kl}  \right )  h_l   =   0  \, ,
\end{equation}
where the array~$\vec h$ contains the components of the yet unknown vectors~$\vec g_l$ sorted (by choice) in the following way
\begin{equation}
\setcounter{MaxMatrixCols}{12}
\vec{h} = 
\begin{pmatrix}
g_1^1 & g_2^1 & \dotsc & g_{q_g+1}^1 & 
g_1^2 & g_2^2 & \dotsc & g_{q_g+1}^2 &
g_1^3 & g_2^3 & \dotsc & g_{q_g+1}^3
\end{pmatrix}^\trans \, .
\end{equation}
It is required that~\eqref{eq:moveLineDesc} is always satisfied irrespective of~$\theta$, which is the case for the right null vectors of the matrix~$C_{kl}$. The right null vectors are determined with a SVD, see e.g.~\cite{strangLinAlg}, yielding the set of null vectors~$\vec g_l^{(i)} $, where the index~$(i)$ denotes the number of the null vector. 

The number of null vectors of~$C_{kl}$ depends on the degrees~$q_x$ and~$q_g$ of the basis~$P_j(\theta)$ and~$\widetilde{P}_l(\theta)$, and the coefficients~$\vec \alpha_j$ of the specific curve considered. For subsequent computations the number of null vectors must be more than the number of intersections of the curve with a line (or its algebraic degree). The non-square matrix~$C_{kl}$ has~$q_x+q_g+1$ rows and~$3(q_g+1)$ columns. Hence, its number of right null vectors must be equal or greater than~$3(q_g+1)-(q_x+q_g+1) = 2 q_g - q_x+2$.\footnote{
%Depending on the coefficients~$\vec \alpha_j$ some of the rows of~$C_{kl}$ might be linearly dependent so that the number of right null vectors is more than~$2 q_g - q_x+2$.
The number of right null vectors is larger when, for instance,  a quadratic curve is described with  a cubic polynomial  (i.e. $q_x=3$ but $\alpha_4=0$). 
More precisely,  if $n$ is the largest integer such that $\vec \alpha_n \neq 0$ then the number of right null vectors of~$C_{kl}$ is exactly~$2q_g- n +3$.
}  
In order to obtain the~$q_x$ intersections it is necessary to have~
\begin{equation}
	 2 q_g - q_x+2 \ge q_x \quad \Rightarrow \quad q_g \ge q_x -1 \, . 
 \end{equation} 
The set of null vectors denoted with~$ \vec g_l^{(i)} $ introduced in~\eqref{eq:lines} yields a set of moving lines 
\begin{equation}~\label{eq:moveLinesFinal}
	l^{(i)} (\theta, \, \vec x) = 
  \begin{pmatrix}
	\vec x  \\ 1
	\end{pmatrix}  
	\cdot 
	\left (  \sum_{l=1}^{q_g+1}  
	\widetilde{P}_l  (\theta) \vec g_l^{(i)}	  \right )
	 =  0
\end{equation}
with a common intersection point on the curve~$\vec x(\theta)$. As an example, in Figure~\ref{fig:movingLines2} the description of a quadratic curve with~$q_x=2$ by two moving lines~$l^{(1)} (\theta, \, \vec x)$  and~$l^{(2)} (\theta, \, \vec x)$ with~$q_g=1$ is shown.

Next, the intersection of a given parametric line 
\begin{equation}
	\vec r(\xi) = \begin{pmatrix} r^1 (\xi) \\ r^ 2 (\xi) \end{pmatrix}= \begin{pmatrix} c_1^1  \\ c_1^2\end{pmatrix} \xi + \begin{pmatrix} c_0^1 \\ c_0^2 \end{pmatrix}
\end{equation}
with the curve~$\vec x(\theta)$ is considered, where $\vec c_1= (c_1^1, \, c_1^2 )^\trans $ and $\vec c_0= (c_0^1, \, c_0^2 )^\trans $ are two given vectors. We require that the determined set of moving lines~$\{ l^{(i)} (\theta, \, \vec x) = 0 \}$ and~$\vec r(\xi)$ intersect at the same point to be yet determined. According to~\eqref{eq:linesDisc}, at the common intersection point it is required that
\begin{equation}
l^{(i)} (\theta, \, \vec r (\xi)) = 
  \begin{pmatrix}
	\vec r(\xi)  \\ 1
	\end{pmatrix}  
	\cdot 
	\left (  \sum_{l=1}^{q_g+1}  
	\widetilde{P}_l  (\theta) \vec g_l^{(i)}	  \right )
	 =  0 \, .
\end{equation}
These equations describe  the intersection of each moving line~$l^{(i)}(\theta, \, \vec x)=0$ with the given line~$\vec r(\xi)$ and can be rewritten as
\begin{equation} \label{eq:intersects-1}
	\sum_{l=1}^{q_g+1} 
	\widetilde{P}_l  (\theta)
	\left [ 
  \begin{pmatrix}
	\vec r(\xi)  \\ 1
	\end{pmatrix}  
	\cdot 
	\left (   \vec g_l^{(i)}	  \right ) \right ]
	 =  0 \, .
\end{equation}
As discussed the number of moving lines satisfies~$\max (i) \ge q_x$ and all of them can be combined in one homogenous equation system   
\begin{equation} \label{eq:intersects}
	\sum_{l=1}^{q_g+1} 
	\widetilde{P}_l  (\theta)
	\left [ 
	 A_{li} - \xi B_{li}  \right ]
	 =  0 \, .
\end{equation}
with~$A_{li}$ and~$ B_{li}$ representing the components of two matrices~$ \vec A$ and~$\vec B $. Both matrices have~\mbox{$q_g+1$} rows and~\mbox{$2q_g-q_x+2$} or more columns. To obtain~$\widetilde{P}_l  (\theta)$ and~$\xi$  that satisfy~\eqref{eq:intersects} the following generalised eigenvalue problem is considered   
\begin{equation}~\label{eq:eigval}
	\vec \phi  \left ( \vec A - \xi \vec B \right ) = \vec 0 \, .
\end{equation} 
That is, the eigenvalues~$ \xi^{(j)}$ are the parametric coordinates of the intersection points on the line~$\vec r(\xi)$ and the eigenvectors are (up to a multiplicative constant) the basis functions~$\widetilde P_l(\theta^{(j)})$ evaluated at the intersection points~$\theta^{(j)}$.
Unfortunately, the matrices~$\vec A$ and~$\vec B$ are not always square and computing~the values $ \xi^{(j)}$ that satisfy this equation requires non-standard linear algebra techniques. However, as will be  discussed in  Section~\ref{sec:example}, for the purposes of intersection computation it is sufficient to consider a square eigenvalue problem obtained from~\eqref{eq:eigval} by taking only some of its columns.  The non-complex eigenvalues of this square eigenvalue problem contain all the intersection points between the given line and the curve. Some of these non-complex eigenvalues may not be actual intersection points, but they can easily be identified. 

A non-complex eigenvalue $\xi^{(j)}$ of~\eqref{eq:eigval}, or its corresponding square  eigenvalue problem, gives the potential intersection point~$\vec r(\xi^{(j)})$. The respective unknown parameter value $\theta^{(j)}$ on the curve satisfies the equation~\mbox{$\vec x(\theta^{(j)}) = \vec r (\xi^{(j)})$}, which is a hard to solve nonlinear problem. According to~\eqref{eq:intersects} and~\eqref{eq:eigval}, however, the left null vector~$\vec \phi^{(j)}$ corresponding~$\xi^{(j)}$  is proportional to the vector~$\widetilde{P}_l  (\theta^{(j)})$, which is exploited to determine~$\theta^{(j)}$.  More precisely, if there is a single parameter value $\theta^{(j)}$, the ratio of any two consecutive components yields
\begin{equation} \label{eq:paramtericEigvec}
	\theta^{(j)} = \frac{\widetilde P_{i+1}}{ \widetilde P_i} = \frac{\phi_{i+1}^{(j)}}{\phi_i^{(j)}} \, .
\end{equation}	
It is assumed here that the monomials in~$\widetilde P_i$ are sorted in increasing order so that the ratio of two consecutive entries is simply~$\theta^{(j)}$. If there are~$p$ parameter values $( \theta^{(j,1)},\, \dotsc,\theta^{(j,p)})$ corresponding to the eigenvalue $\xi^{(j)}$, then  the corresponding left null vector space is given by a matrix $\vec K$ with $p$ rows. Keeping the previous assumption on the ordering of the monomials in~$\widetilde P_i$, we define the \mbox{$p\times p$} matrices $\vec \Delta_i$  by taking the columns $i$ to $p+i$ of $\vec K$. Then, the parameter values $(\theta^{(j,1)}, \, \dotsc,\theta^{(j,p)})$ are obtained by solving  a generalised eigenvalue problem 
\begin{equation}\label{eq:preimagesDelta}
( \vec \Delta_{i+1}-\theta \vec \Delta_i ) \vec  \psi  = \vec 0.	
\end{equation}

%
%--------------------------------------------------------------------------------
\subsection{Illustrative example} \label{sec:example}
\begin{figure}
\centering
	\subfloat[Intersection points given by~\eqref{eq:eigenvalues1}]
	{
		\includegraphics{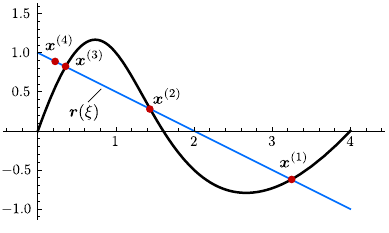}
		\label{fig:curveIntersect1}
	}
	\hfill
	\subfloat[Intersection points given by~\eqref{eq:eigenvalues2}]
	{
		\includegraphics{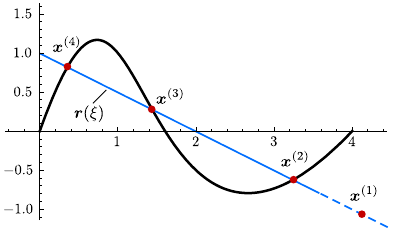}
		\label{fig:curveIntersect2}
	}
\caption{Intersection points between a cubic Lagrange curve and a line.}
\label{fig:curveIntersection}
\end{figure}

%--------------------------------------------------------------------------------
%
The intersection of a cubic Lagrange curve with a line depicted in Figure~\ref{fig:curveIntersection} is considered next. The cubic curve~$\vec x(\theta)$ interpolating the points \mbox{$\vec{x}_1 = (0, \, 0)^\trans$,} \mbox{$\vec{x}_2 = (1, \, 1)^\trans$,} \mbox{$\vec{x}_3 = (2, \, -0.5)^\trans$} and \mbox{$\vec{x}_4 = (4, \, 0)^\trans$} is expressed in power basis with 
%
%\begin{equation}
%	\vec{x}(\theta) = \sum_{j=1}^{4} P_j(\theta) \vec \alpha_j = 
%	\begin{pmatrix}
%	P_1(\theta) \\  P_2(\theta) \\  P_3(\theta) \\ P_4(\theta) 
%	\end{pmatrix}
%	\cdot
%	\begin{pmatrix}
%	\vec{x}_1 \\
%	-5.5\vec{x}_1 + 9\vec{x}_2 -4.5\vec{x}_3 + \vec{x}_4 \\
%	9\vec{x}_1 -22.5\vec{x}_2 + 18\vec{x}_3 - 4.5\vec{x}_4 \\
%	-4.5\vec{x}_1 + 13.5\vec{x}_2 - 13.5\vec{x}_3 + 4.5\vec{x}_4
%	\end{pmatrix}  \, .
%\end{equation}
%
\begin{equation}
	\vec{x}(\theta) = \sum_{j=1}^{4} P_j(\theta) \vec \alpha_j = 
	\begin{pmatrix}
	0 & 4 &  -4.5  &  4.5 \\
	0 & 11.25  & -31.5  & 20.25 \\ 
	\end{pmatrix}  \begin{pmatrix}
	P_1(\theta) \\  P_2(\theta) \\  P_3(\theta) \\ P_4(\theta) 
	\end{pmatrix}\, .
\end{equation}
The intersections between the given line 
\begin{equation} \label{eq:lineEquation}
\vec{r}(\xi) =
\begin{pmatrix}
0  \\  1
\end{pmatrix}
+ \xi
\begin{pmatrix}
4 \\  -2
\end{pmatrix}
=
\begin{pmatrix}
4\xi \\ 1 - 2\xi
\end{pmatrix} 
\end{equation}
and the curve~$\vec x(\theta)$ are sought. 

As discussed, the degree $q_g$ of the auxiliary polynomial needs to satisfy $q_g \geq q_x - 1 = 2$. In this example choosing~$q_g=2$ yields matrices $\vec{A}$ and $\vec{B}$ that have the dimensions $3\times3$ and it is straightforward to compute the eigenvalues~$\xi$ of the generalised eigenvalue problem~\eqref{eq:eigval}. To illustrate the more challenging case with non-square matrices~$\vec{A}$ and~$\vec{B}$, which turns out to be inevitable in the case of surfaces, we choose here~$q_g = 3$   such that 
\begin{equation}
 \vec{g}(\theta) = \sum_{l=1}^{4} \widetilde{P}_l (\theta) \vec g_l  \, .
\end{equation}
The right null vectors~$\vec g_l^{(i)}$ that define the moving lines are obtained from~\eqref{eq:moveLineDesc} with~$\vec C$ and the components
{\normalsize
\begin{equation}
\vec{C} = 
\setcounter{MaxMatrixCols}{12}
\begin{pmatrix*}[r]
0 & 0 & 0 & 0 & 0 & 0 & 0 & 0 & 1 & 0 & 0 & 0 \\
4 & 0 & 0 & 0 & 11.25 & 0 & 0 & 0 & 0 & 1 & 0 & 0 \\
-4.5 & 4 & 0 & 0 & -31.5 & 11.25 & 0 & 
  0 & 0 & 0 & 1 & 0 \\
4.5 & -4.5 & 4 & 0 & 20.25 & -31.5 & 11.25 & 0 & 
  0 & 0 & 0 & 1 \\
0 & 4.5 & -4.5 & 4 & 0 & 20.25 & -31.5 & 11.25 & 0 & 
  0 & 0 & 0 \\
0 & 0 & 4.5 & -4.5 & 0 & 0 & 20.25 & -31.5 & 0 & 0 & 
  0 & 0 \\
0 & 0 & 0 & 4.5 & 0 & 0 & 0 & 20.25 & 0 & 0 & 0 & 0
\end{pmatrix*} \, .
\end{equation} 
}
This matrix has the dimensions~$7 \times 12$, its rank is seven and the  dimension of its right null space is five. Each of the corresponding five independent null vectors $\vec{g}_l^{(i)}$ defines one moving line.
Figure~\ref{fig:curveMovingLines} depicts the first four moving lines $l^{(1)}(\theta, \, \vec x)=0$, $l^{(2)}(\theta, \, \vec x)=0$, $l^{(3)}(\theta, \, \vec x)=0$ and $l^{(4)}(\theta, \, \vec x)=0$ at the parameter values  $\theta=0.4$ and $\theta=0.8$.

Substituting the line~\eqref{eq:lineEquation} as in~\eqref{eq:intersects} gives the two non-square matrices  (with four significant digits)
{\normalsize
\begin{equation} \label{eq:curveMatrixA}
\vec{A} = 
\begin{pmatrix*}[r]
0.08710 & 0.08740 & 0.04298 & 0.05767 & 0.01708 \\
-0.2222 & -0.02814 & 0.8999 & -0.01814 & 0.03355 \\
0.04655 & -0.07068 & -0.03786 & 0.9989 & 0.01874 \\
-0.03786 & 0.1289 & -0.01752 & -0.001102 & 1.010
\end{pmatrix*}
\end{equation}
}
and
{\normalsize
\begin{equation} \label{eq:curveMatrixB}
\vec{B} = 
\begin{pmatrix*}[r]
1.046 & 1.101 & 1.499 & 0.7322 & 0.2201 \\
-3.351 & -1.347 & -0.1816 & 0.04373 & 0.1881 \\
2.390 & -2.402 & -0.1053 & 0.2525 & 0.3777 \\
-0.1302 & 2.987 & -0.3424 & 0.1006 & 0.3987
\end{pmatrix*} \, . 
\end{equation}
}

One approach to obtaining the generalised eigenvalues of matrices $\vec{A}$ and $\vec{B}$ is to use pencil reduction, see \cite{xiaoXiao:2018, BL10, BD}, which is not a widely used linear algebra operation and may introduce additional numerical issues because of several numerical rank estimations. Alternatively, the eigenvalue problems defined by square submatrices $\vec{A}^{\Box}$ and $\vec{B}^{\Box}$ of the largest size, e.g.~the first four columns of~$\vec A$ and~$\vec B$, can be considered
\begin{equation}\label{eq:subEigproblem}
	\vec \phi^{\Box}  \left ( \vec A^{\Box} - \xi \vec B^{\Box} \right ) = \vec 0  \, .
\end{equation}
Although there can only be three intersection points for a cubic curve, this problem has four eigenvalues and eigenvectors. Notice that each column $i$ of~\mbox{$\vec A-\xi \vec B$,} i.e. \mbox{$\sum_l \phi_l (A_{li} - \xi B_{li}) =0$,}   represents the intersection between the line~$\vec r(\xi)$ and the moving line~$l^{(i)}(\theta, \, \vec x) =0$. Taking a different couple of submatrices $\vec{A}^{\Box}$ and $\vec{B}^{\Box}$ of largest size, like the last four columns of~$\vec A$ and~$\vec B$, yields a different set of intersection points. Both sets of intersection points contain the three true intersection points in addition to one fictitious intersection point.

To demonstrate this observation, the intersection points with two different pairs of  largest square submatrices $\vec A^{\Box}$ and $\vec B^{\Box}$ are computed. The eigenvalues for the problem defined by the first four columns of~$\vec A$ and~$\vec B$ are  
\begin{equation} \label{eq:eigenvalues1}
\xi^{(1)} = 0.8112 \, , \quad \xi^{(2)} = 0.3594\, , \quad \xi^{(3)} =0.08875\, , \quad
\xi^{(4)} = 0.05326 \, ,
\end{equation}
and the eigenvalues for the problem defined by the last four columns~$\vec A$ and~$\vec B$ are
\begin{equation} \label{eq:eigenvalues2}
\xi^{(1)} = 28.05\, , \quad \xi^{(2)}  = 0.8112\, , \quad \xi^{(3)} = 0.3594\, , \quad
\xi^{(4)} = 0.08875 \, .
\end{equation}
It is evident that $\xi^{(4)}$ for the first problem and $\xi^{(1)}$ for the second problem correspond to fictitious intersection points while the other three eigenvalues correspond to the true intersection points, see Figure~\ref{fig:curveIntersection}.  The coordinates of the three true intersection points are computed by introducing~$\xi^{(j)}$ in the line equation~\eqref{eq:lineEquation} yielding 
\begin{equation}
	\vec{x}^{(1)} = (3.245, \, -0.6225)^\trans \, , \quad \vec{x}^{(2)} = (1.438, \,  0.2813)^\trans \, , \quad \vec{x}^{(3)} = (0.3550, \,  0.8225)^\trans \, .
\end{equation}
%
%In Figure~\ref{fig:curveIntersect1} and~\ref{fig:curveIntersect2} all the points given in~\eqref{eq:eigenvalues1} and~\eqref{eq:eigenvalues2} are shown. 

The fictitious points can also be detected without computing several eigenvalue problems and comparing their eigenvalues; that is, different from above, by computing only one single eigenvalue problem. This is accomplished by determining the parametric coordinates of the intersection points $\theta^{(j)}$ on the curve~$\vec x(\theta)$. As indicated in~\eqref{eq:paramtericEigvec}, these parametric coordinates are computed using the eigenvectors. As an example, consider the eigenvalue problem~\eqref{eq:subEigproblem} defined by first four columns of~$\vec A$ and~$\vec B$. Its eigenvalues are given in~\eqref{eq:eigenvalues1} and the coordinates of the intersection points are 
\begin{equation} \label{eq:curvePoints}
\begin{aligned}
\vec{x}^{(1)} &= (3.245, \, -0.6225)^\trans , \quad \vec{x}^{(2)} = (1.438, \, 0.2813)^\trans  , \quad
\vec{x}^{(3)}= (0.3550, \, 0.8225)^\trans , \quad  \\ \vec{x}^{(4)} &=
(0.2130, \, 0.8935)^\trans \, .
\end{aligned}
\end{equation}
The respective parameters $\theta^{(j)}$ of the intersection points are
\begin{equation}
\theta^{(1)} = 0.9014\, , \quad \theta^{(2)} = 0.5\, , \quad \theta^{(3)} = 0.09861 \, , \quad \theta^{(4)} = 0.03932 \, .
\end{equation}
 A point is a true intersection point if and only if $\vec{x}(\theta^{(j)}) =  \vec{x}^{(j)}$; otherwise, it is a fictitious point. It can easily be found $\vec{x}(\theta^{(4)}) = (0.1506, \, 0.3949) \neq \vec{x}^{(4)}$ such that the fourth point is not an intersection point. Thus, in practice the intersection points are computed from a single couple of square matrices $\vec{A}^{\Box}$ and $\vec{B}^{\Box}$.

%!TEX root = /Users/lbuse/Inria/Article_Tex/nitro/paper/nonIterative.tex
%
%--------------------------------------------------------------------------------
\section{Intersection of lines with surfaces \label{sec:surfaces}}
%--------------------------------------------------------------------------------
%

The extension of the introduced method to surfaces is straightforward. Let \mbox{$\vec{x} (\vec \theta) = (x^1 (\vec \theta), \, x^2 (\vec \theta), \, x^3(\vec \theta) )^\trans$} be a parametric surface, with~$\vec x(\vec \theta) \in \mathbb R^3$, of bi-degree $(q_x^1, \,  q_x^2)$ given either in Lagrange basis $L_{\vec{i}}(\theta^1, \,\theta^2)$ or power basis~$P_{\vec{j}}(\theta^1, \,  \theta^2)$ with
\begin{equation}
	\vec{x}(\vec{\theta}) = \sum_{\vec{i}}L_{\vec{i}}(\vec{\theta})\vec{x}_{\vec{i}} = \sum_{\vec{j}}P_{\vec{j}}(\vec{\theta})\vec{\alpha}_{\vec{j}} \, ,
\end{equation}
where $\vec{i} = (i^1, \, i^2)$ and $\vec{j} = (j^1, \, j^2)$ are  multi-indices, $\vec{\theta} = (\theta^1, \, \theta^2) $ are the parametric surface coordinates, and $\vec{x}_{\vec{i}} \in \mathbb R^3$ and $\vec{\alpha}_{\vec{j}} \in \mathbb R^3$ are the coefficients in the two basis. Usually, in finite element applications the degrees~$q_x^1$ and~$q_x^2$ of the surface~$\vec x(\vec \theta)$ are the same.

In line with the curve case, a point $\vec{x}(\vec{\theta})$ on the surface is defined as the intersection of several moving planes of the form
\begin{equation} \label{eq:movingPlane}
	l(\vec \theta, \, \vec x) = 
	\begin{pmatrix}
	\vec{x} \\ 1
	\end{pmatrix}
	\cdot \vec{g}(\vec{\theta})
	= x^1g^1(\vec{\theta}) + x^2g^2(\vec{\theta}) + x^3g^3(\vec{\theta}) + g^4(\vec{\theta}) = 0 \, ,
\end{equation}
where $\vec{g}(\vec{\theta}) = \left(g^1(\vec{\theta}), \,  g^2(\vec{\theta}), \, g^3(\vec{\theta}),  \, g^4(\vec{\theta})\right)^\trans$ is an auxiliary vector collecting the parameters of the plane.  Although~$l(\vec \theta, \, \vec x) = 0$ is now a plane instead of a line, it is still denoted with the same symbol to keep the notation simple. The parameters~$\vec g(\vec \theta)$ are assumed to be of the following form 
\begin{equation}
\vec{g}(\vec{\theta}) =
	\sum_{\vec{l}}\widetilde{P}_{\vec{l}}(\vec{\theta})\vec{g}_{\vec{l}} \, . 
\end{equation}
The  bi-degree~$(q_g^1, \, q_g^2)$ of the power basis~$\widetilde{P}(\vec \theta)$ has to be sufficiently high to describe all the intersection points of the surface~$\vec x(\vec \theta)$ with a line. 

The planes describing the surface~$\vec x(\vec \theta)$ have to satisfy   
\begin{equation}
	l(\vec \theta, \, \vec x (\vec \theta)) = 
	\begin{pmatrix}
	\vec{x}(\vec{\theta}) \\[2pt] 1
	\end{pmatrix}
	\cdot
	\vec{g}(\vec{\theta}) = \left[ \sum_{\vec{j}}P_{\vec{j}}(\vec{\theta})
	\begin{pmatrix}\vec{\alpha}_{\vec{j}} \\[2pt]
	1
	\end{pmatrix}\right]
	\cdot
	\left( \sum_{\vec{l}}\widetilde{P}_{\vec{l}}(\vec{\theta})\vec{g}_{\vec{l}}
	\right) = 0 \, ,
\end{equation}
which can  be rearranged to
\begin{equation}
	\sum_{\vec{l}}\sum_{\vec{k}}\widehat{P}_{\vec{k}}(\vec{\theta})C_{\vec{k}\vec{l}}h_{\vec{l}} = 0
\end{equation}
with a new power basis $\widehat{P}_{\vec{k}}$ of bi-degree~$(q_x^1 + q_g^1+1, \, q_x^2 + q_g^2 +1)$ and the array $\vec{h}$ containing the sorted components of the vectors $\vec{g}_{\vec{l}}$. The matrix $C_{\vec{k}\vec{l}}$ has \mbox{$(q_x^1 + q_g^1 + 1)(q_x^2 + q_g^2 + 1)$} rows and \mbox{$4(q_g^1 + 1)(q_g^2 + 1)$} columns. It has at least $4(q_g^1 + 1)(q_g^2 + 1) - (q_x^1 + q_g^1 + 1)(q_x^2 + q_g^2 + 1)$ right null vectors, i.e. the difference between the number of columns and rows.  The surface~$\vec x(\vec \theta)$ has the algebraic degree~$2q_x^1q_x^2$, which is equal to its number of intersections with a line.  In order to obtain all the intersections the following condition has to be satisfied
\begin{equation} \label{eq:minDeg2D}
	4(q_g^1 + 1)(q_g^2 + 1) - (q_x^1 + q_g^1 + 1)(q_x^2 + q_g^2 + 1) \geq 2q_x^1q_x^2 \, .
\end{equation}	
As in the curve case,  the degree of~$\widetilde P (\vec \theta)$ along the $\theta^2$  direction  can be chosen as~$q_g^2 \geq q_x^2 - 1$ yielding 
\begin{equation}
	q_g^1 \geq 2q_x^1 - 1.
\end{equation}
By symmetry, it is also valid to choose $q_g^1 \geq q_x^1 - 1$ and $q_g^2 \geq 2q_x^2 - 1$. After the right null vectors of~$C_{\vec{k}\vec{l}}$ are computed the subsequent steps in computing the intersections are identical to the curve case.

\medskip

Finally, for the sake of completeness, we mention that the introduced intersection algorithm also applies to triangular finite elements. If $\vec x(\vec \theta)$ is a triangular parametric surface of degree $q_x$ then the degree $q_g$ of the auxiliary vector $\vec{g}(\vec{\theta})$ has to be chosen to satisfy $q_g\geq 2(q_x-1)$; see also \cite[\S 3]{laurent2014implicit}.

%
%--------------------------------------------------------------------------------
\section{Conclusions \label{sec:conclusions}}
%--------------------------------------------------------------------------------
%
The introduced method is able to determine in one-shot, simultaneously without iterating, all the intersection points between a line and a curve or surface.  It can be applied to curves or surfaces given in any polynomial basis, like the Lagrange or Bernstein, after a straightforward conversion to the power basis. Aspects requiring further research include cases for which the eigenvalue problem~\eqref{eq:preimagesDelta} is degenerate and the preconditioning of the eigenvalue problem~\eqref{eq:subEigproblem}. In the accompanying implementation, several choices have been considered for both.

\bibliographystyle{wileyj}
\bibliography{nonIterative}

\end{document}